\documentclass{article}
\usepackage{latexsym}
\usepackage{amssymb}
\usepackage{amsfonts}                     
\usepackage{amsmath}
\usepackage[all]{xy}
\usepackage{color}
\usepackage{mathrsfs}
\usepackage[utf8]{inputenc}
\usepackage[T1]{fontenc}
\usepackage{amscd}
\usepackage{amsxtra} 
\usepackage{url}

\newtheorem{thm}{\bf Theorem}[section]
\newtheorem{cor}[thm]{\bf Corollary}
\newtheorem{lem}[thm]{\bf Lemma}
\newtheorem{prop}[thm]{\bf Proposition}
\newtheorem{defn}[thm]{\bf Definition}

\newtheorem{rems}[thm]{\bf Remarks}
\newtheorem{exmp}[thm]{\bf Example}

\newcommand{\field}[1]{\mathbb{#1}}

\newcommand{\Q }{\field{Q}}
\newcommand{\Z }{\field{Z}}

\def\E{{\cal E}}

\def\X{{\cal X}}

\def\Ext{{\rm Ext}}
\def\Hom{{\rm Hom}}

\def\proof{{\parindent0pt {\bf Proof.\ }}}
\def\Tor{{\rm Tor}}
\def\coker{{\rm coker}}

\def\SFP{{\rm $S$-FP}}
\def\FP{{\rm FP}}

\newcommand{\cqfd}
{\hspace{1cm}
\rule{2mm}{2mm}%
\medbreak%
\par%
}

\begin{document}

\title{\SFP-injective modules}

\author{Driss Bennis and Ayoub Bouziri}

\date{}
 
\maketitle
\begin{abstract} Let $R$ be a commutative ring, and let $S$ be a multiplicative subset of $R$. In this paper, we introduce and investigate the notion of \SFP-injective modules. Among other results, we show that, under certain conditions, a ring $R$ is $S$-Noetherian if and only if every \SFP-injective $R$-module is $S$-injective.  Moreover, we establish, under certain conditions, counterparts of Matlis, Stenström and Cheatham-Stone's characterizations of $S$-coherent rings. 
 \end{abstract}

\medskip 
 
{\scriptsize \textbf{Mathematics Subject Classification (2020)}: 13C13, 13E15, 13E99. }
 
 {\scriptsize \textbf{Key Words}: $S$-Noetherian rings,  $S$-coherent rings,  \SFP-injective modules.}


\section{Introduction}
\hskip .5cm  Throughout this paper, $R$ is a commutative ring with identity, all modules are unitary and $S$ is a multiplicative subset of $R$; that is, $1 \in S$ and $s_1s_2 \in S$ for any $s_1 \in S, s_2 \in S$. Unless explicitly stated otherwise, when we refer to a multiplicative subset $S$ of $R$, we implicitly assume that $0 \notin S$. This assumption will be used in the sequel without explicit mention. Let $M$ be an $R$-module. As usual, we use $M^+$ and $M_S$ to denote, respectively, the character module $\Hom_{\Z}(M, \Q/\Z)$ and the localization of $M$ at $S$. Recall that $M_S \cong M \otimes_R R_S$.

 
 \medskip 
 In 1940, Baer initially introduced the notion of injective modules, greatly enriching the study of ring theory. An $R$-module $M$ is said to be injective if $\Ext_R^1(N,M)=0$ for any $R$-module $N$. In 1967, Maddox \cite{Mad1} extended the notion of injective modules to that of absolutely pure modules. Recall that an $R$-module $M$ is said to be absolutely pure if it is a pure submodule of every module containing it. A well-known characterization states that an $R$-module $M$ is absolutely pure if and only if $\Ext^1_R(P, M) = 0$ for any finitely presented $R$-module $P$ \cite[Proposition 2.6]{Ste1}. This characterization led to the terminology of \FP-injective modules instead of absolutely pure modules (\FP \, for finitely presented), which was first used by Stenström \cite{Ste1}.

 \medskip
    Like injective modules, \FP-injective  modules are an important tool for characterizing some classical rings. Megibben \cite{Meg1} showed that a ring $R$ is semihereditary if and only if any homomorphic image of an \FP-injective $R$-module is \FP-injective, and it is Noetherian if and only if any \FP-injective $R$-module is injective. He also proved that a ring $R$ is von Neumann regular if and only if any $R$-module is \FP-injective \cite[Theorem 5]{Meg1}. In 1970, Stenström showed that a ring $R$ is coherent if and only if any direct limit of \FP-injective modules is \FP-injective \cite[Theorem 3.2]{Ste1}. In 1981, Cheatham and Stone, in their work \cite[Theorem 1]{Che1}, established characterization for coherent rings using the notion of character modules as bellow:

  \begin{thm}[\cite{Che1}, Theorem 1]\label{int} The following statements are equivalent for a ring $R$:
\begin{enumerate}
\item $R$ is coherent.
\item $M$ is \text{FP}-injective if and only if $M^+$ is a flat $R$-module.
\item $M$ is \text{FP}-injective if and only if $M ^{++}$ is an injective  $R$-module.
\item $M$ is flat if and only if $M^{++}$ is a flat $R$-module.
\end{enumerate}
\end{thm}

 \medskip
 In 1982, Matlis \cite{Mat1} characterized coherent rings as those rings for which the duality homomorphisms are isomorphisms for certain modules:

 \begin{thm}[\cite{Mat1}, Theorem 1]\label{thm-matlis} Let $R$ be a commutative  ring. The following statements are equivalent:
\begin{enumerate}
\item $R$ is coherent. 
\item The induced morphisms $\Tor^n_R(\Hom_R(N,E), M)\to \Hom_R(\Ext_R^n(M,N), E)$ are isomorphisms for all $n\geq 0$,  whenever $E$ is injective and $M$ is finitely presented.        
\item $\Hom_R(N,E)$ is a flat $R$-module for all injective $R$-modules $N$ and  $E$.
\end{enumerate}
\end{thm}


\medskip
 In the last years, the notion of the $S$-property has drawn the attention of several authors. This notion was introduced in 2002 by D. D. Anderson and Dumitrescu, where they defined the concepts of $S$-finite modules and $S$-Noetherian rings \cite{And2}.  Namely, an $R$-module $M$ is said to be $S$-finite module if there exist a finitely generated submodule $N$ of $M$ and $s \in S$ such that $sM \subseteq N$. A commutative ring $R$ is said to be $S$-Noetherian if every ideal of $R$ is $S$-finite.   To extend coherent rings by multiplicative sets, Bennis and El Hajoui \cite{Ben1} introduced the notions of $S$-finitely presented modules and $S$-coherent rings. An $R$-module $M$ is said to be $S$-finitely presented if there exists an exact sequence of $R$-modules $0 \to K\to L\to  M \to 0$, where $L$ is a  finitely generated free $R$-module and $K$ is an $S$-finite R-module. Moreover, a  commutative ring $R$ is called $S$-coherent, if every finitely generated ideal of $R$ is $S$-finitely presented.  Among other results, they showed that $S$-coherent rings have a similar characterization to the classical one given by Chase for coherent rings \cite[Theorem 3.8]{Ben1}.

\medskip
   There are many notions and results of a homological nature that can be generalized from Noetherian (resp., coherent) rings to $S$-Noetherian rings (resp., $S$-coherent rings).  In this context,  Qi et al. \cite{Qi1}, inspired by the classical case, introduced the notion of $S$-flat modules and provided a new characterization of $S$-coherent rings in terms of these modules.    An $R$-module $M$ is said to be $S$-flat if, for any finitely generated ideal $I$ of $R$, the natural homomorphism $(I \otimes_{R} M)_S \to (R \otimes_{R} M)_S$ is a monomorphism \cite[Definition 2.5]{Qi1}; equivalently, $M_{S}$ is a flat $R_{S}$-module \cite[Proposition 2.6]{Qi1}. Notably, any flat $R$-module is $S$-flat. A general framework for $S$-flat modules was developed in the paper \cite{Bou1}.
   
   \medskip 
    Recently, in \cite{Bou2}, the notion of $S$-injective modules was introduced, and it was shown that many results can be extended from Noetherian rings to $S$-Noetherian rings. For example, it was demonstrated that, under certain conditions on $S$, a ring $R$ is $S$-Noetherian if and only if every direct limit of $S$-injective $R$-modules is $S$-injective, if and only if the class of $S$-injective modules is closed under direct sums \cite[Proposition 2.13 and Corollary 2.17]{Bou3}. Additionally, a counterpart of Cheatham-Stone’s characterizations of $S$-Noetherian rings was presented using the notion of character modules \cite[Theorem 2.22]{Bou3}.

\medskip 
  The primary motivation of this work is to introduce and study an $S$-version of \FP-injective modules, extending some results on coherent and Noetherian rings to $S$-coherent rings and $S$-Noetherian rings, respectively. 
 
 \medskip  
  The organization of the paper is as follows: In Section \ref{s:2}, we introduce the notion of \SFP-injective modules, which serves as a weak version of $S$-injective modules and an $S$-version of \FP-injective modules (see Definition \ref{2-defs}). In Theorem \ref{2-thm-char-sfp-inj}, we provide some characterizations of \SFP-injective modules.  We then proceed to explore some properties of these modules.  For example, we demonstrate that the class of \SFP-injective modules is closed under direct summands, direct products, and, under certain conditions, finite direct sums. Furthermore, we establish that, under certain conditions, the class of all $S$-injective modules coincides with the class of \SFP-injective modules if and only if $R$ is $S$-Noetherian (see Proposition \ref{2-prop-s-inj=sfp-inj-s-neo}). This extends the well-known result that a ring is Noetherian if and only if every \FP-injective module is injective \cite[Theorem 3]{Meg1}. In Section \ref{s:3}, we provide new characterizations for $S$-coherent rings in the case where $R_S$ is finitely presented as an $R$-module.   Firstly, we introduce an $S$-version of Matlis' Theorem in Theorem \ref{thm-s-matlis}. Then, we present an $S$-version of Stenström's characterization of coherent rings in Theorem \ref{thm-s-stens}. Finally, we conclude this section with a counterpart of Cheatham and Stone’s characterizations for $S$-coherent rings in Theorem \ref{thm:cheatham_stone_counterpart}.

\section{\SFP-injective modules}\label{s:2}
\hskip .5cm Recall that an exact sequence $0\to L\to M\to N\to 0$ is pure if and only if, for every finitely presented $R$-module $F$, the induced sequence  $0\to \Hom_R(F,L)\to \Hom_R(F,M)\to\Hom_R(F, N)\to 0$ is exact. An $R$-module $M$ is FP-injective if and only if every monomorphism $M\to N$ is pure; that is, the induced exact sequence   $0\to M\to N\to N/M\to 0$ is pure. Now, we introduce $S$-versions of these notions: 

\begin{defn}\label{2-defs}
\begin{enumerate}
\item An exact sequence of $R$-modules $0\to A\to B\to C\to 0$ is said to be $S$-pure if the induced sequence $0\to \Hom_R(P_S,A)\to \Hom_R(P_S,B)\to \Hom_R(P_S,C)\to 0$  is exact for any finitely presented $R$-module $P$.
\item A monomorphism $f : M\to N$ is said to be $S$-pure if the induced exact sequence $0\to M\to N\to \coker(f)\to 0$ is $S$-pure.
\item A submodule $M$ of an $R$-module $N$ is said to be an $S$-pure  submodule if the inclusion $M\to N$ is an $S$-pure monomorphism.
\item An $R$-module $M$ is said to be \SFP-injective (or absolutely $S$-pure) if every monomorphism $M \to N $ is $S$-pure. 
\end{enumerate}
\end{defn}

 Recall from \cite{Pos1} that an $R$-module $C$ is said to be $S$-weakly cotorsion if $\Ext^1_R(R_S, C) = 0$. An $R$-module $F$ is said to be $S$-strongly flat if $\Ext^1_R(F, C) = 0$ for every $S$-weakly cotorsion $R$-module $C$. We will see that any \SFP-injective $R$ is $S$-weakly cotorsion.  However, the converse does not hold in general (see the discusion just after Theorem \ref{2-thm-char-sfp-inj}).

According to \cite[Definition 2.1]{Bou3}, an $R$-module $E$ is said to be $S$-injective if, whenever $i: A\to B$ is a monomorphism and
$h: A_S\to E$ is any morphism of $R$-modules, there exists a morphism of $R$-modules  $g$ making the following diagram commutes:
$$\xymatrix{  & E & \\
0\ar[r]& A_{S} \ar[r]^{i_{S}}\ar[u]^{h} & B_{S} \ar@{-->}[lu]_{g}}$$

Next, in Corollary \ref{2-cor-of-thm-char-sfp-ijn}, we will see that any $S$-weakly cotorsion $S$-injective $R$-module is \SFP-injective. However, an \SFP-injective $R$-module need not be $S$-injective (see Example \ref{2-expl-sfp-inj-not-s-inj}). In Proposition \ref{2-prop-s-inj=sfp-inj-s-neo}, we will show that if $R$ is $S$-Noetherian, then any \SFP-injective module is $S$-injective. Moreover, the converse holds (at least) when $R_S$ is finitely presented and the $S$-torsion in $R$ is bounded; that is, there exists an element $s_0 \in  S$ such that $sr = 0$ for $s \in S$ and $r \in  R$ implies $s_0r = 0$. This definition can be found in \cite{Pos1}.

\begin{rems}\label{2-rem-def}
\begin{enumerate}
\item Let $P$ be an $S$-finitely presented $R$-module. Then $P_S$ is a finitely presented $R_S$-module by \cite[Remark 3.4 (3)]{Ben1}. By \cite[Proposition 3.19]{Pin2}, there exists a finitely presented $R$-module $P'$ such that $P_S\cong P'_S$. Thus, an exact sequence $0 \to K \to L \to M \to 0$ is $S$-pure if and only if, for every $S$-finitely presented $R$-module $P$, the induced sequence $0 \to \Hom_R(P_S,K) \to \Hom_R(P_S, L) \to \Hom_R(P_S,M) \to 0$  is exact.
\item Let $\X$ be a class of modules. Recall from \cite[Section 1]{War1} that an exact sequence $0 \to K \to L \to M \to 0$ of $R$-modules is called $\X$-pure exact if, for any $X \in \X$, the induced morphism $\Hom_R(X, L) \to \Hom_R(X, M)$ is surjective. Thus, an exact sequence $0 \to K \to L \to M \to 0$  is $S$-pure if and only if it is $\FP_S$-pure exact, where $\FP_S=\{P_S| P \,\text{is finitely presented} \}$.
\item The $U(R)$-pure exact sequences (resp., monomorphisms, submodules) are exactly the pure exact sequences (resp., monomorphisms, submodules), where $U(R)$ denotes the set of units in $R$. 
\item The $U(R)$-\FP-injective modules are exactly the $\FP$-injective modules, where $U(R)$ denotes the set of units in $R$. 
\end{enumerate}
\end{rems}

\begin{prop}
\begin{enumerate}
\item Let $0 \to K \to L \to M \to 0$ be an exact sequence of $R$-modules. Then $0 \to K \to L \to M \to 0$ is $S$-pure if and only if the induced sequence $0\to \Hom_R(R_S, K)\to\Hom_R(R_S, L)  \to \Hom_R(R_S,M)\to 0$  is a pure exact sequence.

\item A monomorphism $K \to L$ is $S$-pure if and only if $0 \to \Hom_R(R_S, N) \to \Hom_R(R_S, M)$ is a pure monomorphism.
\end{enumerate}
\end{prop}

\proof 1. This follows immediately from the following natural isomorphism:\begin{center}
$\Hom_R(A, \Hom_R(B, C))\cong \Hom_R(A \otimes_R B, C)$,
\end{center}
for any $R$-modules $A$, $B$, and $C$ \cite[Theorem 2.75]{Rot1}.

2. This follows from (1). \cqfd

The following lemma will be needed:
\begin{lem}\label{2-lem-ext-of-sub-s-s-inje} Let  $M$ be  an $S$-weakly cotorsion $R$-module, then $$\Ext^1_R(N_S,M) \cong \Ext^1_R(N, \Hom_R(R_S,M)).$$ for any $R$-module $N$.
\end{lem}

\proof  Let $0 \rightarrow M \rightarrow E \rightarrow E/M \rightarrow 0$ be an exact sequence with $E$ is an injective $R$-module. Since $M$ is $S$-weakly cotorsion, the induced sequence
$$0 \rightarrow \Hom_R(R_S, M) \rightarrow \Hom_R(R_S, E) \rightarrow\Hom_R(R_S, E/M) \rightarrow 0$$ is exact. Consider the following  commutative diagram:

\begin{center}
\xymatrix@C=0.4cm@R=0.8cm{
    \Hom_R(N_S,E) \ar[d] \ar[r] & \Hom_R(N_S,E/M) \ar[d] \ar[r] & C \ar@{-->}[d]\ar[r] & 0 \\
    \Hom_R(N, \Hom_R(R_S, E)) \ar[r] & \Hom_R(N,\Hom_R(R_S,E/M)) \ar[r] & D \ar[r] & 0
}

\end{center} where $C=\Ext_R^1(N_S,M)$ and $D=\Ext_R^1(N,\Hom_R(R_S,M))$. The upper row is exact because $E$ is injective. The lower row is exact because $\Hom_R(R_S,E)$ is injective. Since, by \cite[Theorem 2.75]{Rot1}, the first two vertical arrows are isomorphisms, $\Ext_R^1(N_S,M) \cong \Ext_R^1 (N, \Hom_R(R_S,M))$. \cqfd

 Recall that an $R$-module $M$ is $S$-injective if and only if $\Hom_R(R_S, M)$ is injective \cite[Proposition 2.7]{Bou3}. Therefore, the next result follows immediately from Lemma \ref{2-lem-ext-of-sub-s-s-inje}.
\begin{cor}\label{2-cor-char-of-s-inj} The following are equivalent for an $S$-weakly cotorsion $R$-module $M$:\begin{enumerate}
\item $M$ is $S$-injective,
\item $\Ext_R^1(N_S,M)=0$ for any $R$-module $N,$
\item $\Ext_R^1(N,M)=0$ for any $R_S$-module $N.$
\end{enumerate} 
\end{cor}

An $S$-weakly cotorsion $R$-module $M$ need not be $S$-injective (see Corollary \ref{prop-S-absu-colocali-absolu} and Example \ref{2-expl-sfp-inj-not-s-inj}). Recall that an $R$-module $C$ is said to be $S$-$h$-reduced if $\Hom_R(R_S, C) = 0$. Thus, every $S$-$h$-reduced module is $S$-injective, but an $S$-$h$-reduced module is not necessarily $S$-weakly cotorsion.  In fact, the $S$-$h$-reduced $S$-weakly cotorsion modules are the so-called $S$-contramodule $R$-modules, which are extensively studied in papers \cite{Pos2, Pos3}. Whether the assumption in Corollary \ref{2-cor-char-of-s-inj} can be removed remains an interesting question.

Recall that an $R$-module $M$ is \FP-injective if and only if it is a pure submodule of an injective $R$-module. In the next result, we provide a corresponding result for \SFP-injective modules in terms of $S$-pureness and $S$-injective modules. 
\begin{thm}\label{2-thm-char-sfp-inj} The following assertions are equivalent for an $R$-module $M$:
\begin{enumerate}
\item $M$ is \SFP-injective,
\item $M$ is injective with respect to the exact sequences $0 \to A \to B \to C_S \to 0$ with $C$ is  $S$-finitely presented,
\item $M$ is injective with respect to the exact sequences $0 \to A \to B \to C_S \to 0$ with $C$ is finitely presented,
\item $\Ext_R^1(C_S,M) = 0$ for every $S$-finitely presented $R$-module $C$,
\item $\Ext_R^1(C_S,M) = 0$ for every finitely presented $R$-module $C$,
\item $M$ is an $S$-pure submodule of an injective module,
\item $M$ is an $S$-pure submodule of an $S$-weakly cotorsion $S$-injective module.
\end{enumerate}
\end{thm}

\proof 
The equivalences $2. \Leftrightarrow 3.$ and $4. \Leftrightarrow 5.$ follow from the fact that every finitely presented $R_S$-module come form a finitely presented $R$-module \cite[Proposition 3.19]{Pin2}, and the localization $P_S$ of any finitely presented (resp., $S$-finitely presented) $R$-modul $P$  is a finitely presented $R_S$-module.

 $1. \Rightarrow 4.$ Consider an exact sequence $\E = 0 \to M \to N \stackrel{g}{\longrightarrow} C_S \to 0$, where $C$ is an $S$-finitely presented $R$-module. According to $(1)$ and Remarks \ref{2-rem-def}, there exists $h: C_S \to N$ such that $g \circ h = 1_{C_S}$. Therefore, $\E$ splits. Hence, $\Ext_R^1(C_S,M)=0$ by \cite[Theorem 7.31]{Rot1}.

 $2.\Rightarrow 5.$ Similar to $1.\Rightarrow 4.$

 $4. \Rightarrow 1.$ and $5. \Rightarrow 2.$ follow immediately by using the long exact sequences.
 
 $1. \Rightarrow 6.$ This follows from the fact that every $R$-module $M$ can be embedded as a submodule of an injective $R$-module.

$6. \Rightarrow 7.$ This follows from the fact that every injective $R$-module is both $S$-weakly cotorsion and $S$-injective .

$7. \Rightarrow 5.$ If $\E = 0 \to M \to E \stackrel{g}{\longrightarrow} M' \to 0$ is an $S$-pure monomorphism with $E$ being both $S$-weakly cotorsion and $S$-injective,  then, by  Corollary \ref{2-cor-char-of-s-inj}, for any finitely presented $R$-module $C$, the induced sequence $$ \Hom_R(C_S,E)\stackrel{g*}{\longrightarrow}  \Hom_R(C_S,M') \to \Ext_R^1(C_S,M)\to 0$$ is exact. Since $\E$ is $S$-pure, $g^*$ is an epimorphism. Consequently, $\Ext_R^1(C_S,M) = 0$, for any finitely presented $R$-module $C$.  \cqfd

It follows from Theorem \ref{2-thm-char-sfp-inj}(5) that every \SFP-injective $R$-module is $S$-weakly cotorsion. However, an $S$-weakly cotorsion $R$-module need not be \SFP-injective. Indeed, let $S$ be a multiplicative subset of a Noetherian ring $R$ such that $R_S$ is not a semisimple ring.  If every $S$-weakly cotorsion  $R$-module  is \SFP-injective, then for any ideal $I$ of $R$, we have $$\Ext_R^1(R_S/I_S, C)=0$$ for any $S$-weakly cotorsion $R$-module $C$. Hence, $R_S/I_S$ is an $S$-strongly flat $R$-module for any ideal $I$ of $R$. Hence, by \cite[Lemma 3.1]{Baz1}, $R_S/I_S$ is a projective $R_S$-module for any ideal $I$ of $R$. Then, $R_S$ is semisimple, which contradicts our assumption. Thus, there exists an $S$-weakly cotorsion $R$-module that is not \SFP-injective.

\begin{cor}\label{prop-S-absu-colocali-absolu} The following assertion are equivalent for an $R$-module $M$:

\begin{enumerate}
\item $M$ is \SFP-injective, 
\item $M$ is $S$-weakly cotorsion and, $\Hom_R(R_S,M)$ is an FP-injective $R$-module.
\end{enumerate}  \end{cor}
\proof This follows from Lemma \ref{2-lem-ext-of-sub-s-s-inje} and Theorem \ref{2-thm-char-sfp-inj} (5). \cqfd

The following lemma will be needed:
\begin{lem}\label{2-lem-rs-fp-to-fp-as-r-mod} If $R_S$ is finitely presented as an $R$-module, then for every $S$-finitely presented $R$-module $N$, $N_S$ is finitely presented as an $R$-module. 
\end{lem}

\proof Let $0\to K\to F\to N \to 0$ be an exact sequence such that $F$ is a finitely generated free $R$-module and $K$ is $S$-finite; that is, there exists a finitely generated submodule $K'$ of $K$ such that $sK\subseteq K'$ for some $s\in S$. Notice that $K'_S=K_S$. Since $R_S$ is a flat $R$-module, the induced sequence $0\to K_S=K'_S\to F_S\to N_S \to 0$	is exact. Since $R_S$ is finitely generated as $R$-module and  $K_S=K'_S$ is  finitely generated  as an $R_S$-modules, $K_S$ is finitely generated as an $R$-module. Since $R_S$ is finitely presented, then $F_S$, being a direct sum of finitely many copies of $R_S$, is finitely presented as an $R$-module. Therefore, $N_S$ is finitely presented as an $R$-module by \cite[Theorem 2.1.2]{Gla1}.\cqfd

\begin{cor}\label{2-cor-of-thm-char-sfp-ijn}
\begin{enumerate}
\item Assume that $R_S$ is finitely presented as $R$-module, then every \FP-injective $R$-module is  \SFP-injective. 
\item Let $M$ be an $S$-weakly cotorsion $R$-module. If $M$ is $S$-injective, then it is \SFP-injective. 
\end{enumerate}
\end{cor}

\proof 1. Let $M$ be an \FP-injective $R$-module and $N$ a finitely presented $R$-module. Then, $N_S$ is a finitely presented $R_S$-module. Since $R_S$ is finitely presented, by Lemma \ref{2-lem-rs-fp-to-fp-as-r-mod}, $N_S$ is finitely presented as an $R$-module. Then, $\Ext_R^1(N_S, M)=0$. Therefore, $M$ is \SFP-injective by Theorem \ref{2-thm-char-sfp-inj}.

2.     This is clear by Theorem \ref{2-thm-char-sfp-inj} (5) and Corollary \ref{2-cor-char-of-s-inj}.   \cqfd

In  \cite[Proposition 1.2]{Cou1}, Couchot proved that, if $R$ is a coherent ring, then each localization of any \FP-injective $R$-module is an \FP-injective $R_S$-module. Here, we prove  that  the localization $M_S$ of any \SFP-injective $R$-module $M$ is \FP-injective  as an $R_S$-module.  

\begin{prop}\label{2-prop-sfp-inj-to-loc-is-fp-rs-inj} Let $M$ be an $R$-module. If $M$ is \SFP-injective, then $M_S$ is an 
 \FP-injective $R_S$-module.
\end{prop}
\proof Consider a monomorphism $M_S \to E$ where $E$ is an injective $R_S$-module. Let $P$ be a finitely presented $R_S$-module. By \cite[Proposition 3.19]{Pin2}, there exists a finitely presented $R$-module $P'$ such that $P'_S \cong P$.  Then, by using the fact that $\Hom_{R}(A, B)= \Hom_{R_S}(A, B)$  whenever $A$ and $B$ are $R_S$-modules, one can see that $$\Hom_{R_S}(P, E) \to \Hom_{R_S}(P, E/M_S)\to 0$$ is an exact sequence. Thus, $M_S$ is an \FP-injective $R_S$-module. \cqfd

It follows form Proposition \ref{2-prop-sfp-inj-to-loc-is-fp-rs-inj} that, for  a multiplicative subset $S$ of $R$, if every \FP-injective $R$-module is \SFP-injective, then the localization at $S$ of any \FP-injective $R$-module  is an \FP-injective $R_S$-module. Hence, by Corollary \ref{2-cor-of-thm-char-sfp-ijn},  if $R_S$ is finitely presented as an $R$-module, then the localization at $S$ of any \FP-injective  $R$-module  is an \FP-injective $R_S$-module. Moreover, in Proposition \ref{2-prop-rs-fp-inj=s-fp-inj}, we prove that these concepts are equivalent for $R_S$-modules.

\begin{prop}\label{2-prop-rs-fp-inj=s-fp-inj} The following assertions are equivalent for an $R_{S}$-module $M$: 
\begin{enumerate}
\item $M$ is \SFP-injective as an $R$-module.
 \item $M$ is \FP-injective as an $R_S$-module.
\item $M$ is \FP-injective as an $R$-module. 
\end{enumerate} 
\end{prop}

\proof $1. \Rightarrow 2.$ By \cite[Corollary 4.79]{Rot1}, $M$ is naturally isomorphic to its localization $M_S$ as $R_S$-modules. Thus, $M$ is \FP-injective as an $R_S$-module by Proposition \ref{2-prop-sfp-inj-to-loc-is-fp-rs-inj}.

$2. \Rightarrow 3.$ This is \cite[Theorem 3.20]{Pin2}.

$3. \Rightarrow 1.$ Let $M \to E$ be a monomorphism where $E$ is an injective $R_S$-module. Then, by (3), for every finitely presented $R$-module $P$, the morphism $\Hom_R(P,E) \to \Hom_R(P, E/M)$ is an epimorphism. The induced morphism $$\Hom_R(P,E)_S \to \Hom_R(P, E/M)_S$$ is also an epimorphism. By \cite[Lemma 4.87 and Corollary 4.79 (ii)]{Rot1}, $\Hom_R(P,M)_S \cong \Hom_{R_S}(P_S, M)$. Moreover, $\Hom_R(A,B) =\Hom_{R_S}(A,B)$ whenever $A$ and $B$ are $R_S$-modules. Thus, $\Hom_R(P_S,E) \to \Hom_R(P_S, E/M)$ is an epimorphism, and therefore, $M$ is \SFP-injective by Theorem \ref{2-thm-char-sfp-inj} (6).  \cqfd

Megibben \cite{Meg1} gave an equivalent way to define \FP-injective modules:
 An $R$-module $M$ is \FP-injective if and only if every diagram
 $$\xymatrix{  & M & \\
0\ar[r]& P' \ar[r]\ar[u] & P \ar@{-->}[lu]}$$
with $P'$ finitely generated and $P$ projective can be completed to a commutative diagram. Here, we extend this result to \SFP-injective modules.

\begin{thm}   The following statements are equivalent for an $R$-module $M$:
\begin{enumerate}
    \item $M$ is \SFP-injective.
    \item $M$ is $S$-weakly cotorsion, and every diagram
    \[
    \xymatrix{
    & M & \\
    0 \ar[r] & P'_S \ar[r] \ar[u] & P_S \ar@{-->}[lu]
    }
    \]
where $P$ is a finitely generated projective $R$-module and $P'$ is an $S$-finitely generated submodule of $P$, can be completed to a commutative diagram.
    
\item  $M$ is $S$-weakly cotorsion, and every diagram
    \[    \xymatrix{    & M & \\    0 \ar[r] & P'_S \ar[r] \ar[u] & P_S \ar@{-->}[lu]
    }
    \]

where $P$ is a finitely generated projective (free) $R$-module and $P'$ is a finitely generated submodule of $P$, can be completed to a commutative diagram.
\end{enumerate}
\end{thm}

\proof $1. \Rightarrow 2.$ Assume that $M$ is \SFP-injective. The induced exact sequence:
\begin{center}
$0\to P'_S\to P_S \to (P'/P)_S\to 0$.
\end{center}
gives rise to the exact sequence:
\begin{center}
$ \Hom(P_S, M) \to \Hom(P'_S, M)\to Ext_R^1((P'/P)_S,M)$.
\end{center}
Since $M$ is \SFP-injective and $P/P'$ is $S$-finitely presented, $\text{Ext}_R^1((P/P)_S,M)=0.$ Thus the sequence $$ \Hom(P_S, M) \to \Hom(P'_S, M)\to 0$$ is exact, and therefore the appropriate diagram can be completed.

$2. \Rightarrow 3.$ Trivial.

$3.\Rightarrow 1.$  We show that $\Hom_R(R_S,M)$ is \text{FP}-injective.  By Corollary \ref{prop-S-absu-colocali-absolu}, it follows that $M$ is \SFP-injective. Let $P'$ be a finitely generated submodule of a finitely generated free $R$-module $F$. Then we have the exact sequence $0\to P\to F\to F/P\to 0$ which gives rise to the commutative diagram with exact rows:
  \begin{center}
  
$\xymatrix{
  \Hom_R(F_S,M)\ar[d] \ar[r]&\Hom_R(P_S,M)   \ar[d]\ar[r]  & 0 \\
\Hom_R(F, \Hom_R(R_S, M))\ar[r]& \Hom_R(P,\Hom_R(R_S, M)) &}  $
   \end{center} 
   The upper row is exact by assumption. Moreover, all the vertical maps are isomorphisms. Hence, $\Hom_R(F, \Hom_R(R_S, M)) \to \Hom_R(P,\Hom_R(R_S,M))$ is an epimorphism. Therefore, $\Hom_R(R_S,M)$ is \text{FP}-injective as desired.\cqfd

 According to \cite[Theorem 4 and Remark on page 565]{Meg1}, if $R$ is a coherent ring, then an $R$-module $M$ is \FP-injective if and only if $\Ext_R^1(R/I, M) = 0$ for any finitely generated ideal $I$ of $R$. Replacing "coherent" with "$S$-coherent," we prove the following theorem, which will be useful in the next Section.
 
Recall from \cite[Definition 3.3]{Ben1} that a ring $R$ is called $S$-coherent, if every finitely generated ideal of $R$ is $S$-finitely presented.
 
 \begin{thm}\label{2-th-s-version-mjeben}
Suppose that $R$ is $S$-coherent. Then the following assertion  are equivalent for an $R$-module $M$:

 \begin{enumerate}
  \item $M$ is \SFP-injective.
  \item $M$ is $S$-weakly cotorsion, and  $\Ext_R^1(R_S/I_S,M) =0$ for all $S$-finitely generated ideals $I$ of $R$.
  \item $M$ is $S$-weakly cotorsion, and $\Ext_R^1(R_S/I_S,M) =0$ for all finitely generated ideals $I$ of $R$.
\end{enumerate}
 \end{thm}
 \proof $1.\Rightarrow 2.$ and $2.\Rightarrow 3.$ are straightforward.
 
 $3.\Rightarrow 1.$ Assume that $\Ext_R^1(R_S/I_S, M) = 0$ for all $S$-finitely generated ideals $I$. Consider an arbitrary homomorphism $f:N_S \to M$, where $N$ is a finitely generated submodule of a finitely generated free $R$-module $F$.  We need to show that $f$ can be extended to $F_S$.   By Theorem \ref{2-th-s-version-mjeben}, if follows that $M$ is \SFP-injective. 
 We  prove this by induction on the number $n$ of generators of $F$. The case $n =1$ holds by assumption. For $n \geq 2$, we write $F = Rx \oplus V$, where $V$ is free with $(n-1)$ generators.

Let $I= \{r\in R/ rx\in N+V\}$. The morphism $rx+g\in N \to r\in I$, induces an isomorphism $I\cong N/(N\cap V)$. Consequently,  $I$ is finitely generated. Since $R$ is $S$-coherent, $I$ is $S$-finitely presented. Hence, by \cite[Theorem 2.4]{Ben1}, $(N\cap V)$ is $S$-finite; that is, there exists a finitely generated  submodule $K$ of $N\cap V$   such that $s(N\cap V)\subseteq K$ for some $s\in S$. In particular, $(N\cap V)_S=K_S$.  Thus, we can assume that $N\cap V$ is finitely generated.

 By induction, there exists a homomorphism  $g : V_S\to M$ such that $g\mid_{(N\cap V)_S} = f\mid_{(N\cap V)_S}$. Hence, there is  a unique homomorphism $h: (N+V)_S\to M$ extending both $g$ and $f$. Define $\theta: I_S\to M$ by $\theta(\frac{r}{s})=h(\frac{rx}{s})$ for all $r\in I$. Since $\Ext_R^1(R_S/I_S,M) =0$, there is a homomorphism $\psi : R_S\to M$ that extends $\theta$. The mapping $k : F_S\to M$ defined by $k(\frac{rx+v}{s})=\psi(\frac{r}{s})+g(\frac{v}{s})$ for all $r\in R$, $s\in S$, and  $v\in V$ is a homomorphism extending $f$.  \cqfd

There are questions that arise naturally about classes of modules. Two of these are whether the direct sum of elements of the class remains in the class and whether the direct limit of elements in the class also remains in the class. The direct sum of \FP-injective modules is always \FP-injective, but it is not easy to see whether this holds for \SFP-injective modules. Here, we examine specific cases, but we leave the general case as an open question.\begin{prop}\label{2-prop-prod-of-s-inj}
Let $(M_i)_{i\in I}$ be a family of $R$-modules. Then
\begin{enumerate}
\item  $\prod\limits_{i\in I} M_i$ is \SFP-injective if and only if each $M_i$ is \SFP-injective. In particular, every direct summand of an \SFP-injective module is \SFP-injective.

\item  Assume that $R_S$ is finitely presented as an $R$-module. Then $\bigoplus\limits_{i\in I} M_i$ is \SFP-injective if and only if each $M_i$ is \SFP-injective. 

\end{enumerate}
\end{prop}

\proof 1. 
By Theorem \ref{2-thm-char-sfp-inj}, $\prod\limits_{i\in I} M_i$ is \SFP-injective if and only if, for any finitely presented $R$-module $P$, $\Ext_R^1(P_S, \prod\limits_{i\in I} M_i)=0$. By \cite[Proposition 7.22]{Rot1}, we have $\Ext_R^1(P_S, \prod\limits_{i\in I} M_i) \cong \prod\limits_{i\in I}\Ext_R^1(P_S, M_i)$. Thus, $\Ext_R^1(P_S, \prod\limits_{i\in I} M_i)=0$ if and only if, for each $i\in I$, $\Ext_R^1(P_S, M_i)=0$. Again, by Theorem \ref{2-thm-char-sfp-inj}, this holds if and only if, for each $i\in I$, $M_i$ is \SFP-injective.

2. Since $R_S$ is flat, it is projective \cite[Theorem 2.6.15]{Wan1}. Hence, every $R$-module is $S$-weakly cotorsion. Thus, by Corollary \ref{prop-S-absu-colocali-absolu}, $\bigoplus\limits_{i\in I} M_i$ is \SFP-injective if and only if $\Hom_R(R_S, \bigoplus\limits_{i\in I} M_i)$ is FP-injective. By \cite[Theorem 2.6.10]{Wan1}, we have $$\Hom_R(R_S, \bigoplus\limits_{i\in I} M_i)\cong \bigoplus\limits_{i\in I}\Hom_R(R_S, M_i).$$ Since the class of FP-injective modules is closed under direct sums, $\bigoplus\limits_{i\in I}\Hom_R(R_S, M_i)$ is FP-injective if and only if, for each $i\in I$, $\Hom_R(R_S, M_i)$ is FP-injective. Again, by Corollary \ref{prop-S-absu-colocali-absolu}, this happens if and only if, for each $i\in I$, $M_i$ is \SFP-injective.  \cqfd

Megibben \cite[Theorem 3]{Meg1} showed that a ring $R$ is Noetherian if and only if every \text{FP}-injective $R$-module is injective. Here, we extend this result in some particular cases. The general case is left as an open question.

\begin{prop}\label{2-prop-s-inj=sfp-inj-s-neo}  If $R$ is $S$-Noetherian, then every \SFP-injective module is $S$-injective. Moreover, if the $S$-torsion in $R$ is bounded and $R_S$ is finitely presented as an $R$-module, then the converse holds true.
\end{prop}

\proof  
Assume that $R$ is $S$-Noetherian. Let $M$ be an \SFP-injective $R$-module. Let $I$ be an ideal of $R$. Since $R$ is $S$-Noetherian, $I$ is $S$-finite.  Hence, $\Ext_R^1(R_S/I_S,M)=0$. Therefore, $M$ is $S$-injective by \cite[Proposition 2.4]{Bou3}.

Conversely, suppose that every \SFP-injective module is $S$-injective. If $R$ is not $S$-Noetherian, then by \cite[Corollary 2.16]{Bou3}, there exists a family $(M_i)_{i\in I}$ of $S$-injective modules such that $\bigoplus\limits_{i\in I} M_i$ is not $S$-injective. However, by Proposition \ref{2-prop-prod-of-s-inj}, $\bigoplus\limits_{i\in I} M_i$ is \SFP-injective, a contradiction. Thus, $R$ is $S$-Noetherian. \cqfd

Now, utilizing Proposition \ref{2-prop-s-inj=sfp-inj-s-neo}, we can demonstrate that the classes of \SFP-injective modules and $S$-injective modules are not coincident in general.
 
\begin{exmp}\label{2-expl-sfp-inj-not-s-inj}
Let $R'$ be a commutative ring and $S'$ be a multiplicative subset of $R'$ such that the $S'$-torsion in $R'$ is bounded, and $R'_{S'}$ is finitely presented (for an example of such a ring, please see \cite[Example 2.23]{Bou3}). Let $R''$ be a non-Noetherian commutative ring. Consider the ring $R = R'\times R''$ with the multiplicative subset $S = S' \times U(R'')$, where $U(R'')$ denotes the set of units in $R''$. Then, there exists an \SFP-injective $R$-module which is not $S$-injective.
\end{exmp}

\proof  Firstly, notice that the $S$-torsion in $R$ is bounded, and $R_S = (R')_{S'} \times R''$ is a finitely presented $R$-module. Moreover, since $R''$ is not ($U(R'')$-)Noetherian,  $R$ is not $S$-Noetherian. Hence, according to Proposition \ref{2-prop-s-inj=sfp-inj-s-neo}, there exists an \SFP-injective module which is not $S$-injective. \cqfd

\section{S-coherent rings}\label{s:3} 
\hskip .5cm  In this section, we provide some characterizations of $S$-coherence. Specifically, considering particular multiplicative subsets $S$, we present counterparts  of Matlis, Stenström, and Cheatham-Stone’s characterizations of $S$-coherent rings.

         Recall that an $R$-module $M$ is said to be $S$-flat if for any finitely generated ideal $I$ of $R$, the natural homomorphism $I \otimes_{R} M \to R \otimes_{R} M$ is an $S$-monomorphism; equivalently, $M_{S}$ is a flat $R_{S}$-module \cite[Proposition 2.6]{Qi1}.

We need the following lemmas:

\begin{lem}\label{s-flat-of-prod-and-quoti}
Let $(F_i)_{i\in I}$ be a family of $S$-flat modules. For each $i  \in I$, let $K_i$ be  a submodule of $F_i$ such that $F_i/K_i$ is $S$-flat. Suppose that $R_S$ is finitely presented. Then,the following assertions are equivalent: \begin{enumerate}
\item $\prod\limits_{i\in I}F_i$ is $S$-flat. 
\item $\prod\limits_{i\in I}K_i$ and $\prod\limits_{i\in I}F_i/K_i$ are $S$-flat. 
\end{enumerate}
\end{lem}
\proof
 $2. \Rightarrow 1.$ Let $N$ be a finitely presented  $R$-module.  Consider the induced exact sequences: $$\Tor_R^1(\prod\limits_{i\in I}K_i, N_S)\to \Tor_R^1(\prod\limits_{i\in I}F_i, N_S)\to Tor_R^1(\prod\limits_{i\in I}F_i/K_i, N_S)$$

By \cite[Proposition 2.5]{Bou2}, the first and last terms are zero, which implies that $\Tor_R^1(\prod\limits_{i\in I}F_i, N_S)=0$. Hence, $\prod\limits_{i\in I}F_i$ is $S$-flat, again by  \cite[Proposition 2.5]{Bou2}.

$1. \Rightarrow 2.$  Let $N$ be a finitely presented  $R$-module. Consider the following commutative diagram:
$$\xymatrix@C=0.4cm@R=0.8cm{ 
0 \ar[r] & \Tor_R^1(\prod\limits_{i\in I}F_i/K_i, N_S) \ar[r] \ar@{-->}[d]^{} &  (\prod\limits_{i\in I}K_i)\otimes_R N_S \ar[d]^{\cong} \ar[r] &(\prod\limits_{i\in I}F_i )\otimes_R N_S \ar[d]^{\cong}  \\
 & \prod\limits_{i\in I}\Tor_R^1(F_i/K_i, N_S)   \ar[r] & \prod\limits_{i\in I}(K_i\otimes_R N_S) \ar[r] &   \prod\limits_{i\in I}(F_i\otimes_R N_S)   }$$
 
Since $R_S$ is finitely presented,  $N_S$ is a  finitely presented $R$-module by Lemma \ref{2-lem-rs-fp-to-fp-as-r-mod}. Hence, by \cite[Theorem 3.2.22]{Eno2},  the right two vertical homomorphisms are isomorphisms.  Since,  for each $i\in I$, $F_i/K_i$ is $S$-flat, we have $\Tor_R^1(F_i/K_i, N_S)=0$. Hence,  $\prod\limits_{i\in I}\Tor_R^1(F_i/K_i, N_S)=0$, and so $\Tor_R^1(\prod\limits_{i\in I}F_i/K_i, N_S)=0$. Therefore, $\prod\limits_{i\in I}F_i/K_i$ is  $S$-flat by  \cite[Proposition 2.5]{Bou2}. 

Now, using \cite[Proposition 2.5]{Bou2} and the induced exact sequence
$$\Tor_R^2(\prod\limits_{i\in I}F_i/K_i, N_S)\to \Tor_R^1(\prod\limits_{i\in I}K_i, N_S)\to \Tor_R^1(\prod\limits_{i\in I}F_i, N_S)$$ we deduce that $\prod\limits_{i\in I}K_i$  is $S$-flat too.  \cqfd 

Recall that if $ E$ is an injective cogenerator $R$-module (universal injective in the terminology of \cite{Meg1}), then  $R \subseteq \Hom_R(E,E)$ and $\Hom_R(E,E)$ is flat if and only if $\Hom_R(E,E)/R $ is flat. To provide an $S$-version of Matlis's characterization for coherent rings, we need an $S$-version of this result:

\begin{lem}\label{len-H-iff-R/H} Let $E$ be an  be an injective cogenerator $R$-module and let $H = \Hom_R(E,E)$ then, $H$
is $S$-flat if and only if $H/R$ is $S$-flat.
\end{lem}
\proof The "if" part follows from the fact that the class of $S$-flat modules is closed under extension \cite[Lemma 2.1]{Bou2}. For the "only if" part, let $I_S$ be an ideal  of $R_S$, where  $I$ is an ideal of $R$. By \cite[Lemma 2.6.5]{Gla1}, $HI\cap R=I$. Hence, $(HI\cap R)_S=H_SI_S\cap R_S=I_S$. Then,  by \cite[Theorem 1.2.3]{Gla1}, $H_S$ is a flat $R_S$-module. \cqfd

\begin{thm}\label{thm-s-matlis} Let $R$ be a commutative ring and $S \subseteq R$ be a multiplicative subset such that $R_S$ is finitely presented as an $R$-module. Then the following assertions are equivalent: 
\begin{enumerate}
\item $R$ is $S$-coherent.
\item $R_S$ is a coherent ring.
\item The induced morphisms $\Tor^n_R(\Hom_R(N,E), M_S)\to \Hom_R(\Ext_R^n(M_S,N), E)$ are isomorphisms for all $n\geq 0$ whenever $E$ is injective and $M$ is $S$-finitely presented. 
       
\item The induced morphisms $\Tor^n_R(\Hom_R(N,E), M_S)\to \Hom_R(\Ext_R^n(M_S,N), E)$ are isomorphisms for all $n\geq 0$ whenever $E$ is injective and $M$ is finitely presented.  
\item $\Hom_R(N,E)$ is an $S$-flat $R$-module for all $S$-\text{FP}-injective $R$-modules $N$ and  any injective $R$-module  $E$. 
   
\item $\Hom_R(N,E)$ is an $S$-flat $R$-module for all $S$-injective $R$-modules $N$ and  any injective module  $E$. 
\item $\Hom_R(N,E)$ is an $S$-flat $R$-module for all injective $R$-modules $N$ and any injective module  $E$.
\end{enumerate}
\end{thm}
\proof  $1.\Rightarrow 2.$ This is \cite[Remarks 3.4(3)]{Ben1}.

$2.\Rightarrow 3.$ Let $M$ be an $S$-finitely presented $R$-module. By \cite[Remark 3.4]{Ben1}, $M_S$ is a finitely presented $R_S$-module.  Since $R_S$ is a coherent ring, $M_S$ has a projective resolution composed of finitely generated $R_S$-modules \cite[Corollary 2.5.2]{Gla1}. 

On the other hand, as $R_S$ is a finitely generated projective $R$-module \cite[Theorem 3.56]{Rot1}, every finitely generated projective $R_S$-module is also a finitely generated projective $R$-module. Therefore, $M_S$ has a projective resolution composed of finitely generated $R$-modules. Consequently, (3) follows from \cite[Theorem 1.1.8]{Gla1}. 
  
  $3.\Rightarrow 4.$ and $6.\Rightarrow 7.$ are obvious. 
  
  
  $4.\Rightarrow 5.$ This follows from Theorem \ref{2-thm-char-sfp-inj} and \cite[Proposition 2.5]{Bou2}. 

  $5.\Rightarrow 6.$ This follows from Corollary \ref{2-cor-of-thm-char-sfp-ijn}.

 $7. \Rightarrow 1.$ Let $E$ be an injective cogenerator. Then, $H := \Hom_R(E, E)$ is $S$-flat by (7). Hence, by Lemma \ref{len-H-iff-R/H},  $H/R$ is $S$-flat. 
 
  Let $I$ be an index set. Denote, for each $i\in I$, $R_i=R$, $H_i=H$, and $E_i=E$. We have, by (4),
  $\prod\limits_{i\in I} H_i=\Hom_R(E, \prod\limits_{i\in I} E_i)$ is $S$-flat because $\prod\limits_{i\in I} E_i$ is injective.  Then $\prod\limits_{i\in I} R_i$ is $S$-flat by Lemma \ref{s-flat-of-prod-and-quoti} and therefore, $R$ is $S$-coherent by \cite[Theorem 4.4]{Qi1}. \cqfd

Next, we present an $S$-version of  Stenström's characterizations of coherent rings.

\begin{thm}\label{thm-s-stens} Suppose that $R_S$ is  finitely presented. Then the following assertions are equivalent: 
\begin{enumerate}
\item $R$ is $S$-coherent.
\item Every direct limit of \SFP-injective  $R$-modules is an \SFP-injective $R$-module. 
\item Every direct limit of  \FP-injective  $R_S$-modules is an \FP-injective $R_S$-module. 
\end{enumerate}
\end{thm}

\proof $1.\Rightarrow 2.$ Let $(M_i)_{i\in J}$  be a direct system of \SFP-injective modules over a directed set $J$.  Let $I$ be a finitely generated ideal of $R$.
We construct the following commutative diagram:
 $$\xymatrix@C=0.4cm@R=0.8cm{ 
\varinjlim \Hom_R( R_S, M_i) \ar[r]\ar[d] &  \varinjlim \Hom_R( I_S, M_i)\ar[r] \ar[d]^{} & \varinjlim \Ext_R^1(R_S/I_S, M_i) \ar@{-->}[d]^{} \ar[r] &0  \\
\Hom_R( R_S, \varinjlim M_i) \ar[r] & \Hom_R( I_S, \varinjlim M_i) \ar[r] & \Ext_R^1(R_S/I_S, \varinjlim M_i) \ar[r] &   0   }$$
 
 Since $R$ is $S$-coherent, $I$ is $S$-finitely presented. Hence, $I_S$ is finitely presented $R_S$-module. Given that $R_S$ is finitely presented, by Lemma \ref{2-lem-rs-fp-to-fp-as-r-mod}, $I_S$ is also a finitely presented $R$-module. Therefore, the first tow vertical morphisms in the above diagram are isomorphisms. Consequently, by the Five Lemma, $\Ext_R^1(R_S/I_S, \varinjlim M_i)\cong \varinjlim \Ext_R^1(R_S/I_S, M_i)=0$. Thus, $\varinjlim M_i$ is \SFP-injective by Theorem \ref{2-th-s-version-mjeben}.

 $2.\Rightarrow 3.$ This follows from Proposition \ref{2-prop-rs-fp-inj=s-fp-inj}.
 
 $3.\Rightarrow 1.$ By \cite[Theorem 3.2]{Ste1}, $R_S$ is a coherent ring. Since $R_S$ is  finitely presented, $R$  is $S$-coherent by Theorem \ref{thm-s-matlis}. \cqfd

Finally, we present an $S$-counterpart to the classical result by Cheatham and Stone [7, Theorem 1]. To do so, recall the following lemma:

\begin{lem}[\cite{Bou3}, Lemma  2.5.2]\label{2-lem-lambek}
An $R$-module $M$ is $S$-flat if and only if its character $\Hom_{\Z}(M,\Q/\Z)$ is $S$-injective.
\end{lem}

\begin{thm}\label{thm:cheatham_stone_counterpart} Let $R$ be a ring and $S$ a multiplicative subset of $R$.   Consider the following assertions:
\begin{enumerate} 
\item $R$ is $S$-coherent,
 \item $M$ is \SFP-injective if and only if $M^+$ is $S$-flat,
 \item $M$ is \SFP-injective if and only if $M ^{++}$ is $S$-injective,
 \item $M$ is $S$-flat if and only if $M^{++}$ is $S$-flat,
\end{enumerate}
The implications $2.\Rightarrow 3.\Rightarrow 4.\Rightarrow 1.$ hold true. Assuming that $R_S$ is finitely presented as an $R$-module, then all the four  assertions are equivalent.

\end{thm}

 \proof 
 $2.\Rightarrow 3.$ This follows from Lemma \ref{2-lem-lambek}.

$3.\Rightarrow 4.$ Let $M$ be an $S$-flat $R$-module. Then, by Lemma \ref{2-lem-lambek},  $M^+$ is $S$-injective; so, by $(3)$ $M ^{+++}$ is $S$-injective and hence, again by Lemma \ref{2-lem-lambek}, $M^{++}$ is $S$-flat. 

Conversely, if $M^{++}$ is $S$-flat, then $M$, being a pure submodule of $M^{++}$ \cite[Proposition 2.3.5]{Xu1}, is $S$-flat by \cite[Lemma 2.1]{Bou1}.

$4.\Rightarrow 1.$   Let $(M_i)_{i\in I}$ be a family of flat $R$-modules. We aim to show that $\prod\limits_{i\in I} M_i$ is $S$-flat. It follows then from \cite[Theorem 4.4]{Qi1} that $R$ is $S$-coherent.  Since every flat $R$-odule is $S$-flat, $\bigoplus\limits_{i\in I} M_i$ is $S$-flat, and hence \begin{center}
$(\bigoplus\limits_{i\in I} M_i)^{++}\cong (\prod_{i\in I} M_i^+)^+$
\end{center} is $S$-flat by (4).  By \cite[Lemma 1]{Che1}, $\bigoplus\limits_{i\in I} M_i^+$ is a pure submodule of $\prod M_i^+$. Then the induced sequence $$(\prod M_i^+)^{+}\to (\bigoplus M_i^+)^+ \to 0$$ splits. Thus, $\prod M_i^{++} \cong (\bigoplus M_i^+)^+$ is $S$-flat. Since, for each  $i\in I$, $M_i$ is a pure submoduele of $M^{++}$, $\prod M_i$ is a pure submodule of $\prod M_i^{++}$ by \cite[Lemma 1]{Che1}. Hence, $\prod\limits_{i\in I} M_i$ is $S$-flat by \cite[Lemma 2.1]{Bou1}.

$1.\Rightarrow 2.$ Assume $R_S$ is finitely presented. Since $R$ is $S$-coherent, $R_S$ is a coherent ring. Then, by \cite[Lemma 2.21]{Bou3},    we have the following isomorphism:
$$ \Tor_R^1(M^+, (R/I)_S) \cong \Ext_R^1((R/I)_S,M)^+$$ for any finitely generated ideal $I$ of $R$.  Thus, (2) follows from Theorem \ref{2-th-s-version-mjeben} and \cite[Proposition 2.5]{Bou2}. \cqfd

Let $R$ be an $S$-perfect ring; that is, every $S$-flat $R$-module is projective \cite[Definition 4.1]{Bou1}. Then, $R_S$ is a finitely presented $R$-module. Indeed, $R_S$ is projective and cyclic as an $R$-module \cite[Theorem 4.9]{Bou1}. Therefore, by \cite[Proposition 3.11]{Rot1}, $R_S$ is finitely presented. Thus, Theorems \ref{thm-s-matlis}, \ref{thm-s-stens} and  \ref{thm:cheatham_stone_counterpart} characterize, in particular, when an $S$-perfect ring $R$ is $S$-coherent. We conclude with the following example:

\begin{exmp}\label{exmp-fin}
Let $R_1$ be a semisimple ring and let $S_1$ be a multiplicative subset of $R_1$.  Let $R_2$ be a non-coherent commutative  ring. Consider the ring $R=R_1\times R_2$ with the multiplicative subset $S = S_1\times 0$. Then, 

\begin{enumerate}
\item $R_S$ is finitely presented.
\item $R$ is  $S$-coherent, but it  not is coherent.
\end{enumerate}  

\end{exmp} 

\proof $1.$ Since $R_1$ is $S_1$-perfect, $(R_1)_{S_1}$ is a finitely generated projective $R_1$-module by \cite[Theorem 4.9]{Bou1}. Then $R_S \cong (R_1)_{S_1}\times 0$ is a finitely generated projective $R$-module, hence it is finitely presented.

$2.$ $R$ is $S$-coherent by \cite[Proposition 3.5]{Ben1}. \cqfd

\noindent\textbf{}

Driss Bennis:  Faculty of Sciences, Mohammed V University in Rabat, Rabat, Morocco.

\noindent e-mail address: driss.bennis@fsr.um5.ac.ma; driss$\_$bennis@hotmail.com

Ayoub Bouziri: Faculty of Sciences, Mohammed V University in Rabat, Rabat, Morocco.

\noindent e-mail address: ayoubbouziri66@gmail.com; ayoub$\_$bouziri@um5.ac.ma

\end{document}